\documentclass{amsart}

\usepackage[english]{babel}
\usepackage{leftidx}
\usepackage[numbers]{natbib}
\usepackage{bigints}
\usepackage{empheq}
\usepackage{array}
\usepackage{graphicx}
\usepackage{lipsum}
\usepackage{amsthm}
\usepackage{needspace}
\usepackage[makeroom]{cancel}
\newenvironment{myproof}[2] {\paragraph{\textbf{Proof of {#1} {#2} :}}}{\hfill$\square$}
\newtheorem{theorem}{Theorem}[section]
\newtheorem{proposition}[theorem]{Proposition}
\newtheorem{lema}[theorem]{Lemma}
\newtheorem{question-non}[]{}
\newtheorem{notation}[theorem]{Notation}
\newtheorem{cor}[theorem]{Corollary}
\newtheorem{definition}[theorem]{Definition}
\newtheorem{observation}[theorem]{Remark}
\newtheorem{example}[theorem]{Example}

\title{On Warped Product Gradient Yamabe Solitons}

\author{Tokura, W. $^{1}$}
\address{$^{1}$ Universidade Federal de Goi\'as, IME, 131, 74001-970, Goi\^ania, GO, Brazil.}
\email{williamisaotokura@hotmail.com $^{1}$}

\author{Adriano, L. $^{2}$}
\address{$^{2}$ Universidade Federal de Goi\'as, IME, 131, 74001-970, Goi\^ania, GO, Brazil.}
\email{levi@ufg.br $^{2}$}

\author{Pina, R. $^{3}$}
\address{$^{3}$ Universidade Federal de Goi\'as, IME, 131, 74001-970, Goi\^ania, GO, Brazil.}
\email{romildo@ufg.br $^{3}$}

\author{Barboza, M. $^{4}$}
\address{$^{4}$ Insituto Federal Goiano, 75790-000, Rodovia Geraldo Silva Nascimento Km 2,5, Uruta\'i, GO, Brazil.}
\email{marcelo.barboza@ifgoiano.edu.br $^{4}$}

\thanks{$^{1,4}$ Supported by CAPES}

\keywords{Warped product, gradient Yamabe solitons, scalar curvature, semi-Riemannian metric, Almost gradient Yamabe solitons.}

\subjclass[2010]{53C21, 53C50, 53C25} 

\begin{document}

\begin{abstract}The purpose of this article is to study gradient Yamabe soliton on warped product manifolds. First, we prove triviality results  in the case of noncompact base with limited warping function, and for compact base. In order to provide nontrivial examples, we consider the base conformal to a semi-Euclidean space, which is invariant under the action of a translation group, and then we characterize steady solitons. We use this method to give infinitely many explicit examples of complete steady gradient Yamabe solitons.
\end{abstract}
\maketitle
\section{Introduction and main results}
\label{intro}

A \textit{Yamabe soliton} is a semi-Riemannian manifold $(M,g)$
admitting a smooth vector field $X\in \mathfrak{X}(M)$ such that
\begin{equation}\label{eq:01}(S_{g}-\rho)g=\frac{1}{2}\mathfrak{L}_{X}g,
\end{equation}
where $S_{g}$ denotes the scalar curvature of $M$, $\rho$ is a real
number and $\mathfrak{L}_{X}g$ denotes the Lie derivative of $g$ with respect to $X$. We say that $(M,g)$ is
\textit{shrinking}, \textit{steady} or \textit{expanding}, if $\rho>0$, $\rho=0$ , $\rho<0$,
respectively. When $X=\nabla h$ for some smooth function $h\in
C^{\infty}(M)$, we say that $(M,g,\nabla h)$ is an
\textit{gradient Yamabe soliton} with \textit{potential function} $h$. In
this case the equation \eqref{eq:01} turns out
\begin{equation}\label{eq8}(S_{g}-\rho)g=Hess_{g}(h),
\end{equation}
where $Hess_{g}(h)$ denote the Hessian of $h$. When $h$ is constant, we
call it a \textit{trivial Yamabe soliton}.

Adding the condition of constant $\rho$ in definition 
\eqref{eq:01} to be a differentiable function on $M$ we obtain the extension of Yamabe solitons called \textit{Almost Yamabe soliton}. In particular, for gradient vector field, we call it a \textit{Almost gradient Yamabe soliton} \cite{barbosa2013conformal}.

After their introduction in the Riemannian sense, the study of semi-Riemannian Yamabe solitons attracted a growing number of authors, showing many differences with respect to the Riemannian case. Calvi\~no \textit{et. al.} in \cite{calvino2012three} study Yamabe solitons and left-invariant Ya\-ma\-be soliton on three-dimensional homogeneous space and showed that the class of Yamabe solitons is strictly greater than the class of left-invariant Yamabe solitons. On the other hand, Neto \textit{et. al.} in \cite{neto2018gradient} showed that the class of semi-Riemannian Yamabe solitons is striclty greater than the Riemannian Yamabe soliton. Moreover, infinitely many proper semi-Riemannian gradient Yamabe solitons are exhibit.



In recent works on solitons and equations, the notion of warped product introduced by Bishop and O'neil in \cite{o1983semi} has attracted major research activities \cite{al2012warped, barboza2018invariant, blaga2017warped, de2017gradient, feitosa2017construction, kim2013warped, leandro2017invariant}.

\begin{definition}(\cite{o1983semi}) Let $(B^{n}, g_{B})$ and $(F^{d}, g_{F})$ be two semi-Riemannian manifolds,  as well as a positive
	smooth function $f$ on $B$. On the product manifold $B\times F$, we define the metric
	\begin{equation}\label{warped metric}
	\overline{g}=\pi^{\ast}g_{B}+(f\circ\pi)^{2}\sigma^{\ast}g_{F},
	\end{equation}
	where $\pi:B\times F\rightarrow B$, $\sigma:B\times F\rightarrow F$ are the projections on the first and second factor, respectively. The product space $B\times F$ furnished with metric tensor $\overline{g}$ is called \textit{warped product}. We denote it by $B\times_{f}F$. The function $f$ is called \textit{warping function}, $B$ is called the \textit{base} and $F$ the \textit{fiber}.
\end{definition}

Brozos-V\'azquez \textit{et al.} in \cite{brozos2016local} provide a special warped product structure for gradient Yamabe solitons, its results establish that a gradient Yamabe soliton $(M,g)$ with potential function $h$ and such that $|\nabla h|\neq0$, is locally isometric to a warped product of unidimensional base and constant scalar curvature fiber. In the Riemannian context a global structure result was given in \cite{cao2012structure}.

Sousa and Pina investigated gradient Ricci soliton on warped product  and proved that the potential function only depends on base or the warping function is constant (see \cite{de2017gradient}). The same technique can be used to prove this result for gradient Yamabe solitons.

These results make it interesting for further investigation of gradient Yamabe solitons with warped product structure $B\times_{f}F$ where the potential function only depends on $B$, and $F$ with constant scalar curvature.

\begin{notation}
	Throughout this paper, we will consider the following:
	\vspace{0,1cm}
	\begin{equation*}
	\overline{M}^{n+d}=\left(B^{n}\times F^{d},\bar{g}\right),\hspace{0,5cm}
	S_{F}=\lambda_{F}=\text{constant}, \hspace{0,5cm} \widetilde{h}=h\circ\pi,\hspace{0,5cm}h\in C^{\infty}(B).
	\end{equation*}
	where $S_{F}$ is the scalar curvature of $F$, $\widetilde{h}$ is the potential function of $\overline{M}^{n+d}$ and $\overline{g}$ is given by \eqref{warped metric}.
\end{notation}

We start by focusing our attention on compact base gradient Yamabe soliton $\overline{M}^{n+d}$. It has been known that every compact Riemannian gradient Yamabe soliton
is of constant scalar curvature, hence, trivial since $h$ is
harmonic, see \cite{daskalopoulos2013classification}, \cite{hsu2012note}. In this direction, we have the following theorem.

\begin{theorem}\label{trivial}
	Let $\overline{M}^{n+d}$ be a gradient Yamabe soliton with compact Riemannian base. Then $\overline{M}^{n+d}$ is trivial.
\end{theorem}

By the above theorem, in order to obtain examples of nontrivial gradient Yamabe solitons, we must relax the hypothesis of compactness of the base. Since the continuous images of compact spaces are compact, that is, limited functions, it then follows that a condition weaker than compactness is considered to be limited. Then we consider the following question.

\textbf{Question:}
Does there exist a gradient Yamabe soliton $\overline{M}^{n+d}$ with  nonconstant limited warping function?

In the sequel, we give a negative partial answer as follows:

\begin{theorem}\label{teorema2}
	Let $\overline{M}^{n+d}$ be a gradient Yamabe soliton with soliton constant $\rho$ and Riemannian base with scalar curvature $S_{B}\geq\rho-\frac{\lambda_{F}}{f^2}$. If $f$ reaches the maximum,  then $\overline{M}^{n+d}$ must be a standard semi-Riemannian product.
\end{theorem}

Next, it is interesting to know under which conditions an gradient Yamabe soliton $\overline{M}^{n+d}$ has non compact base. In this case, we obtain the following

\begin{proposition}\label{proposicao}
	Let $\overline{M}^{n+d}$ be a gradient Yamabe soliton with potential function $\widetilde{h}$ and complete Riemannian base $(B^{n},g_{B})$. If $\langle\nabla \log f,\nabla h \rangle=\text{constant}\neq0$, then $(B^{n},g_{B})$ is isometric to the standard Euclidean space $(\mathbb{R}^{n},g_{0})$.
\end{proposition}

Recently, the conformal semi-Euclidean space has become an interesting space to display examples of steady gradient Yamabe solitons and steady gradient Ricci solitons. Neto and Tenenblat in \cite{neto2018gradient} treated the conformally flat semi-Riemannian space
$(\mathbb{R}^{n},\frac{1}{\varphi^{2}}g_{0})$, where $g_{0}$ is the
canonical semi-Riemannian metric, and obtain a necessary and sufficient
condition to this manifold be a gradient Yamabe soliton. In order to exhibit solutions they consider the invariant action of
an $(n-1)$-dimensional translation group, and as a result geodesically complete(see definition \ref{defi}) example was obtained. The same technique is used to obtain all invariant solutions of steady gradient Ricci soliton \cite{barbosa2014gradient}.

The warped product has proved its efficiency in constructing new examples of manifolds with certain geometric characteristics \cite{ besse2007einstein, dobarro1987scalar, ganchev2000riemannian}. Considering invariant solutions,  Neto in \cite{leandro2017invariant} provide an explicit example of a complete static vacuum Einstein space-time. On the other hand, in the same invariant solution context Sousa in \cite{de2017gradient} provide examples of non conformally flat Ricci solitons.

In the remainder of this article, we focus our attention on the warped product $\overline{M}^{n+d}$ where the base is conformal to an $n$-dimensional semi-Euclidean space, invariant under the action on an $(n-1)$-dimensional translation group. As application, we will construct 5 examples of steady gradient yamabe solitons. Besides, we provide a way to construct infinitelly many explicit examples of geodesically complete steady gradient Yamabe solitons, with base conformal to the Lorentzian space (see example \ref{exemplo completo}).

More precisely, consider the semi-Riemannian metric
\begin{equation*}
\delta=\sum_{i=1}^{n}\varepsilon_{i}dx_{i}\otimes dx_{i}
\end{equation*}
in coordinates $x=(x_{1},\dots,x_{n})$ of $\mathbb{R}^{n}$, where $n\geq 3$, $\varepsilon_{i}=\pm1$. For an arbitrary choice of non zero vector $\alpha=(\alpha_{1},\dots,\alpha_{n})$ we define the function $\xi:\mathbb{R}^{n}\rightarrow\mathbb{R}$ by
\begin{equation*}
\xi(x_{1},\dots,x_{n})=\alpha_{1}x_{1}+\dots+\alpha_{n}x_{n}.
\end{equation*}

Next, we consider that $\mathbb{R}^{n}$ admits a group of symmetries consisting of translations \cite{olver2000applications} and we then look for positive smooth functions $\varphi, f, h:(a,b)\subset\mathbb{R}\rightarrow(0,\infty)$ such that $f=f\circ\xi, \varphi=\varphi\circ\xi, h=h\circ\xi: M=\xi^{-1}(a,b)\rightarrow\mathbb{R}$, satisfies \eqref{eq8} with $M=\mathbb{R}^{n}\times F^{d}$ and metric tensor
\begin{equation*}
g=\frac{\delta}{\varphi(x)^2}+f(x)^2g_{F}.
\end{equation*}

\begin{theorem}\label{teorema invariancia geral}With $(\mathbb{R}^{n},\delta)$ and $f=f\circ\xi, \varphi=\varphi\circ\xi, h=h\circ\xi$ as above, the manifold $M=\mathbb{R}^{n}\times F^{d}$, furnished with the metric tensor
	\begin{equation}\label{metrica com base Rn}
	g=\frac{\delta}{\varphi(x)^2}+f(x)^2g_{F}
	\end{equation}
	is a gradient Yamabe soliton if, and only if, 
	\begin{equation}\label{eq:09}
	\vspace{0,6cm}h''+2\frac{\varphi'h'}{\varphi}=0,
	\end{equation}
	\begin{equation}\label{eq1}
	\begin{split}
	||\alpha||^{2}\Big{[}(n-1)(2\varphi\varphi''-n(\varphi')^{2})-2\frac{d}{f}(\varphi^{2}f''-(n-2)\varphi\varphi'f')-&\frac{d(d-1)}{f^2}\varphi^{2}(f')^{2}+\\&+
	\varphi'h'\varphi\Big{]}=\rho-\frac{\lambda_{F}}{f^{2}},
	\end{split}
	\end{equation}
	\begin{equation}\label{eq2}
	\begin{split}
	||\alpha||^{2}\Big{[}(n-1)(2\varphi\varphi''-n(\varphi')^{2})-2\frac{d}{f}(\varphi^{2}f''-(n-2)\varphi\varphi'f')-&\frac{d(d-1)}{f^2}\varphi^{2}(f')^{2}+\\&-\frac{\varphi^{2}}{f}f'h'\Big{]}=\rho-\frac{\lambda_{F}}{f^{2}},
	\end{split}
	\end{equation}
	when $||\alpha||^{2}\neq0$, that is, $\alpha$ is a timelike or spacelike vector.
	
	And
	
	\begin{equation}\label{eq:10}
	\vspace{0,6cm}h''+2\frac{\varphi'h'}{\varphi}=0,
	\end{equation}
	\begin{equation}\label{eq10}
	\rho-\frac{\lambda_{F}}{f^{2}}=0,
	\end{equation}
	when $||\alpha||^{2}=0$, that is, $\alpha$  is a lightlike vector.
\end{theorem}

\begin{cor}\label{corolario de  rigidez na funcao torcao}If $||\alpha||^{2}=0$ and $\lambda_{F}\neq0$, then the warped product $\mathbb{R}^{n}\times_{f} F^{d}$ become a standard semi-Riemannian product.
\end{cor}

In this case, we have the following obstrution on the constant soliton $\rho$

\begin{cor}If $||\alpha||^{2}=0$ and $\lambda_{F}>0$, then there is no expanding or
	steady gradient Yamabe soliton with product metric \eqref{metrica com base Rn}. Similarly, if we assume that $||\alpha||^{2}=0$ and $\lambda_{F}<0$, then there is no shrinking or
	steady gradient Yamabe soliton with product metric \eqref{metrica com base Rn}.
\end{cor}

Now, by equations \eqref{eq:09} and \eqref{eq:10} in Theorem \ref{teorema invariancia geral} we easily see that a necessary condition for
the above manifold be a gradient
Yamabe soliton, is that
$h$ is a monotone function. That is,
\begin{equation}h'(\xi)=\frac{k_{1}}{\varphi^{2}(\xi)},\nonumber
\end{equation}
for some $k_{1}\in\mathbb{R}$.

We provide steady solutions for ODE in Theorem \ref{teorema invariancia geral} in the cases:
$h'=0$ and $h'\neq0$ with $n+d=6$, or $||\alpha||^{2}=0$.

\begin{theorem}\label{eq5}If $||\alpha||^{2}\neq0$ and $\lambda_{F}\neq0$, then the warped product metric $g$ given by \eqref{metrica com base Rn} with $n+d=6$ is a steady gradient Yamabe soliton with potential function $h$ with $h'\neq0$ if, and only if, the functions $f$, $h$ and $\varphi$, satisfies
	\begin{equation}\label{eqa}f(\xi)=\frac{k_{2}}{\varphi(\xi)},
	\end{equation}
	\begin{equation}\label{eqb}h(\xi)=k_{1}\bigintsss\frac{1}{\varphi^{2}(\xi)}d\xi,
	\end{equation}
	\begin{equation}\label{eqc}\bigintsss\frac{pd\varphi}{q\varphi^{3}W\big(k_{3}e^{\frac{p^2}{4q\varphi^4}-1}\big)+q\varphi^{3}}=\xi+k_{4},
	\end{equation}
where $k_{1}\neq0$, $k_{2}\neq0$, $k_{3}\neq0$, $k_{4}$ are constants, $p=\dfrac{k_{1}}{10}$, $q=\dfrac{\lambda_{F}}{10k_{2}^{2}||\alpha||^{2}}$ and $W$ is the product log function.
\end{theorem}

\begin{theorem}\label{eq5'}If $||\alpha||^{2}\neq0$ and $\lambda_{F}=0$, then the warped product metric $g$ given by \eqref{metrica com base Rn} with $n+d=6$ is a steady gradient Yamabe soliton with potential function $h$ with $h'\neq0$ if, and only if, the functions $f$, $h$ and $\varphi$, satisfies
	\begin{equation}\label{eqa'}f(\xi)=\frac{k_{2}}{\varphi(\xi)},
	\end{equation}
	\begin{equation}\label{eqb'}h(\xi)=k_{1}\bigintsss\frac{1}{\varphi^{2}(\xi)}d\xi,
	\end{equation}
	\begin{equation}\label{eqc'}40\int\frac{\varphi d\varphi}{k_{1}-20k_{3}\varphi^{4}}=\xi+k_{4},
	\end{equation}
	where $k_{1}\neq0$, $k_{2}\neq0$, $k_{3}$ and $k_{4}$, are constants.
\end{theorem}

\begin{theorem}\label{eq7} If $||\alpha||^{2}\neq0$ and $\lambda_{F}=0$, then given a smooth function $\varphi>0$ the warped product metric $g$ given by \eqref{metrica com base Rn} is a steady gradient Yamabe soliton with potential function $h$ with $h'=0$ if, and only if, the functions $h$ and $f$, satisfies
	\begin{equation}\label{eq23}h(\xi)=\text{constant},
	\end{equation}
	\begin{equation}\label{eq22}f(\xi)=\varphi^{\frac{n-2}{d+1}}(\xi)e^{\Phi(\xi)}\left(\int e^{-(d+1)\Phi(\xi)}d\xi+\frac{2}{d+1}C\right)^{\frac{2}{d+1}},
	\end{equation}
where $\Phi(\xi)=\int z_{p}(\xi)d\xi$ and $z_{p}$ is a particular solution of
\begin{equation}\label{Ricati}
z^2+\frac{2}{d+1}z'+\frac{(n+d-1)}{d(d+1)^2}\left(n\left(\frac{\varphi'}{\varphi}\right)^2-2\frac{\varphi''}{\varphi}\right)=0.
\end{equation}
	
\end{theorem}

In the null case $||\alpha||^{2}=0$ we obtain

\begin{theorem}\label{eq6} If $||\alpha||^{2}=0$ and $\lambda_{F}=0$, then given two smooth
	functions $\varphi(\xi)$ and $f(\xi)$, the warped product metric $g$ given by \eqref{metrica com base Rn} is a steady gradient Yamabe soliton with potential function $h$ if, and only if
	\begin{equation}h(\xi)=k_{1}\int\frac{1}{\varphi^{2}(\xi)}d\xi.\nonumber
	\end{equation}
	
\end{theorem}

\begin{observation}As we can see in the proof of Theorem \ref{teorema invariancia geral}, if $\rho$ is a function defined only on the
	base, then we can easily extend Theorem \ref{teorema invariancia geral} into context of
	almost gradient Yamabe solitons. In the particular case of lightlike
	vectors there are infinitely many solutions, that is, given
	$\varphi$ and $f$
	\begin{equation}\rho(\xi)=\frac{\lambda_{F}}{f(\xi)^{2}},\nonumber
	\end{equation}
	\begin{equation}h(\xi)=k_{1}\int\frac{1}{\varphi^{2}(\xi)}d\xi,\nonumber
	\end{equation}
	provide a family of almost gradient Yamabe soliton with warped
	product structure.

\end{observation}

Before proving our main results, we present some examples
illustrating the above theorems.

\begin{example}In Theorem \ref{eq5}, consider $\overline{M}^{6}=\mathbb{R}^{3}\times\mathbb{H}^{3}$, $k_{1}=1$, $k_{2}=1$, $k_{3}=0$, $k_{4}=0$, then we get
	\begin{equation*}f(\xi)=\frac{1}{\varphi(\xi)},\quad h(\xi)=\bigintsss\frac{1}{\varphi^{2}(\xi)}d\xi,\quad \bigintsss\frac{d\varphi}{\varphi^{3}W\big(e^{-\frac{1}{4\varphi^4}-1}\big)+\varphi^{3}}=-10\xi.
	\end{equation*}	
The family of solutions of $\varphi$ is described in the follow phase portrait
	\begin{figure}[!htb]
		\centering
		\includegraphics[scale=0.7]{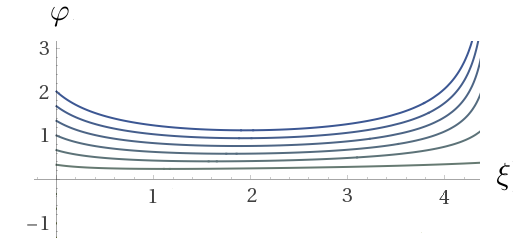}
		\caption{Sampling $\varphi(0)$ and $\varphi'(0)$ }
		\label{Label de referência para a imagem}
	\end{figure}

\end{example}

\begin{example}In Theorem \ref{eq5'}, consider $k_{1}=1$, $k_{2}=1$, $k_{3}=0$, $k_{4}=0$, then the functions
	\begin{equation*}f(\xi)=\sqrt{\frac{20}{\xi}}, \quad h(\xi)=20\ln\xi,\quad \varphi(\xi)=\sqrt{\frac{\xi}{20}},
	\end{equation*}
	provide a steady gradient Yamabe soliton defined in the semi-space $\xi>0$ of Euclidean space $\mathbb{R}^{6-d}$, $d=1,2,3$.

\end{example}

\begin{example}In Theorem \ref{eq5'}, consider $k_{1}=1$, $k_{2}=1$, $k_{3}=-\frac{1}{20}$, $k_{4}=0$, then the functions
	\begin{equation*}f(\xi)=\sqrt{\frac{1}{\tan(\frac{\xi}{20})}}, \quad h(\xi)=20\ln\big{[}\sin\big{(}\frac{\xi}{20}\big{)}\big{]},\quad \varphi(\xi)=\sqrt{\tan\big{(}\frac{\xi}{20}\big{)}},
	\end{equation*}
provide a steady gradient Yamabe soliton defined in the slice $0<\xi<10\pi$ of Euclidean space $\mathbb{R}^{6-d}$, $d=1,2,3$.

\end{example}

\begin{example}In Theorem \ref{eq7} consider the Lorentzian space $(\mathbb{R}^{4},g)$ with coordinates $(x_{1},x_{2},x_{3},x_{4})$
	and signature $\varepsilon_{1}=-1$, $\varepsilon_{k}=1$ for $k=2,3,4$, and
	$F^{3}$ scalar flat fiber. Let $\xi=x_{2}+x_{3}+x_{4}$ and choose
	$\varphi(\xi)=|\sec(\xi)|$ where $\xi\neq\frac{\pi}{2}+k\pi$,
	$k\in\mathbb{Z}$, then $z_{p}(\xi)=-\frac{1}{2}$ is a particular
	solution of $\eqref{Ricati}$, and by Theorem \ref{eq7}
	\begin{equation}f(\xi)=(2|\sec(\xi)|e^{\xi})^{\frac{1}{2}}, \hspace{0,3cm} h(\xi)=\text{constant},\hspace{0,3cm}\varphi(\xi)=|\sec(\xi)|,\nonumber
	\end{equation}
	provide a steady gradient Yamabe soliton in warped metric defined on
	$x_{1}+x_{2}+x_{3}\neq\frac{\pi}{2}+k\pi$, $k\in\mathbb{Z}$.

\end{example}

\begin{example}\label{exemplo completo}In Theorem \ref{eq6} consider the Lorentzian space $(\mathbb{R}^{n},g)$ with coordinates $(x_{1},\dots,x_{n})$
	and signature $\varepsilon_{1}=-1$, $\varepsilon_{i}=1$ for all $i\geq2$, and fiber $(\mathbb{R}^{d},g_{0})$ where $g_{0}$ is the Euclidean metric. Let $\xi=x_{1}+x_{2}$ and
	choose $k\in\mathbb{R}\setminus\{0\}$. Then
	\begin{equation}f(\xi)=e^{k\xi},\hspace{0,3cm}h(\xi)=-\frac{k_{1}e^{-2k\xi}}{2k},\quad k_{1}\neq0,\hspace{0,3cm}\varphi(\xi)=e^{k\xi},\nonumber
	\end{equation}
	defines a family of complete steady gradient Yamabe soliton on
	$(\mathbb{R}^{n},\varphi^{-2}g)\times_{f}(\mathbb{R}^{d},g_{0})$ with potential
	function $h$ and warping function $f$(see section \ref{proofs}).
\end{example}

\section{Preliminaries}
\label{preliminaries}

In this section we shall present some preliminaries which will be used in the paper. We shall follow the notation and terminology of Bishop and O'Neill \cite{o1983semi}.


\begin{lema}\label{lemma}
	$\overline{M}^{n+d}$ is a \textbf{gradient Yamabe soliton} with potential function $\widetilde{h}$ and soliton constant $\rho$ if, and only if, $(B^{n},g_{B})$ is an \textbf{almost gradient Yamabe soliton} with potential function $h$, soliton function 
	\begin{equation*}
	\lambda=-\frac{\lambda_{F}}{f^2}+\frac{2d}{f}\Delta f+d(d-1)\frac{|\nabla f|^2}{f^2}+\rho,
	\end{equation*}
	and scalar curvature
	\begin{equation}\label{escalar}
	S_{B}=\frac{\langle\nabla f,\nabla h\rangle}{f}+\lambda.
	\end{equation}
\end{lema}

\begin{myproof}{Lemma}{\ref{lemma}}Using that $S_{F}=\lambda_{F}$
, we have by the well known formula of scalar curvature on warped product that
\begin{equation*}
S_{\bar{g}}=\pi^{\ast}\left[S_{B}+\frac{\lambda_{F}}{f^2}-2d\frac{\Delta f}{f}+d(d-1)\frac{|\nabla f|^2}{f^2}\right]\\[5pt]
\end{equation*}
where $\Delta$ denote the Laplacian on $B$, and we use $g_{B}=\langle \cdot, \cdot\rangle=|\cdot|^2$ for simplicity.

Then we have that $\overline{M}^{n+d}$ is a gradient Yamabe soliton with potential function $\widetilde{h}$ and soliton constant $\rho$ if, and only if,
	\begin{equation}\label{eq88}\left(\pi^{\ast}\left[S_{B}+\frac{\lambda_{F}}{f^2}-2d\frac{\Delta f}{f}-d(d-1)\frac{|\nabla f|^2}{f^2}\right]-\rho\right)\overline{g}=Hess(\widetilde{h}).\\[5pt]
	\end{equation}

	Let $\mathcal{L}(B)$, $\mathcal{L}(F)$ the spaces of lifts of vector fields on $B$ and $F$ to $B\times F$, respectively. Consider $X\in\mathcal{L}(B)$ and $V\in\mathcal{L}(F)$, then
	\begin{equation*}
	\overline{g}(X,V)=0=Hess(\overline{h})(X,V).
	\end{equation*}
	Hence, we just need to look at equation \eqref{eq88} for pair of fields in $\mathcal{L}(B)$, and $\mathcal{L}(F)$.
	
	Taking $X, Y\in\mathcal{L}(B)$ in \eqref{eq88} and using $Hess(\widetilde{h})=\pi^{\ast}(Hess (h))$ we obtain the following equivalent condition
		\begin{equation*}\pi^{\ast}\left(S_{B}+\frac{\lambda_{F}}{f^2}-2d\frac{\Delta f}{f}-d(d-1)\frac{|\nabla f|^2}{f^2}-\rho\right)\pi^{\ast}g_{B}=\pi^{\ast}(Hess (h)),\\[5pt]
	\end{equation*}
	which says that $(B^{n},g_{B})$ is an almost gradient Yamabe soliton with soliton function
	\begin{equation*}
	\lambda=-\frac{\lambda_{F}}{f^2}+\frac{2d}{f}\Delta f+d(d-1)\frac{|\nabla f|^2}{f^2}+\rho,
	\end{equation*}
and potential function $h$.

Now, consider $V, W\in\mathcal{L}(F)$, then we obtain the Hessian expression
\begin{eqnarray}\label{eq13}
Hess(\widetilde{h})(V,W)& = & V(W(\widetilde{h}))-(\nabla_{V}W)^{M}\widetilde{h}\nonumber\\
&=& V(W(\widetilde{h}))+\frac{\overline{g}( V,
	W)}{\widetilde{f}}\nabla(\widetilde{f})(\widetilde{h})-\nabla_{V}^{F}W(\widetilde{h})\nonumber\\
&=&\widetilde{f}\sigma^{\ast}g_{F}(V,W)(\nabla\widetilde{f})\widetilde{h}\nonumber\\
&=&\widetilde{f}\sigma^{\ast}g_{F}(V,W)[d\pi(\nabla \widetilde{f})(h)\circ\pi].
\end{eqnarray}

Substituting $V, W\in\mathcal{L}(F)$  into \eqref{eq88} and considering equation \eqref{eq13} we obtain
	\begin{equation*}\pi^{\ast}\left(S_{B}+\frac{\lambda_{F}}{f}-2d\frac{\Delta f}{f}-d(d-1)\frac{|\nabla f|^2}{f^2}-\rho\right)\widetilde{f}^2\sigma^{\ast}g_{F}=\widetilde{f}\sigma^{\ast}g_{F}[d\pi(\nabla \widetilde{f})(h)\circ\pi]\\[5pt]
\end{equation*}
And using $\nabla f=\pi_{\ast}(\nabla \widetilde{f})$ we obtain the following equivalent condition
	\begin{equation*}
S_{B}+\frac{\lambda_{F}}{f^2}-\frac{2d}{f}\Delta f-d(d-1)\frac{|\nabla f|^2}{f^2}-\rho= \frac{\langle\nabla f, \nabla h\rangle}{f}
	\end{equation*}
	which is equation \eqref{escalar}. This concludes the proof.
\end{myproof}

\begin{definition}\label{defi}A semi-Riemannian manifold for which every geodesic is defined on the entire real line is said to be \textit{geodesically complete}, or just
\textit{complete}.
\end{definition}

Given a curve $\gamma$ in $M\times_{f}F$, we can write $\gamma(s)=(\gamma_{B}(s),\gamma_{F}(s))$, where $\gamma_{B}=\pi\circ\gamma$ and $\gamma_{F}=\sigma\circ\gamma$. The following proposition guarantees a condition for curve $\gamma$ to be geodesic.

\begin{proposition}\label{geodesicas produto torcido}(\cite{o1983semi}) A curve $\gamma=(\gamma_{B}, \gamma_{F})$ in $B\times_{f}F$ is a geodesic if, and only if,
	\\[0.1pt]
	\begin{enumerate}
		\item $\gamma_{B}''=g_{F}(\gamma_{F}',\gamma_{F}')f\circ\gamma_{B}\nabla f$\quad \text{in B},\\[0.1pt]
		\item $\gamma_{F}''=\dfrac{-2}{f\circ\gamma_{B}}\dfrac{d(f\circ\gamma_{B})}{ds}\gamma_{F}'$ \quad \text{in F}.
	\end{enumerate}
\end{proposition}

\section{Proof of the main result}
\label{proofs}

\begin{myproof}{Theorem}{\ref{trivial}} By Lemma \ref{lemma}, we have that 
	\begin{equation}\label{fundamental}
	(S_{B}-\lambda)g_{B}= Hess (h),\qquad S_{B}-\lambda= \frac{\langle\nabla f, \nabla h\rangle}{f}.
	\end{equation}
	
	Combining this equations we obtain
	\begin{equation}\label{operador}
	\Delta h-\langle\nabla w, \nabla h\rangle=0,
	\end{equation}
	where $w=\ln f^{n}$.
	
		Denoting $\Delta_{w}:=\Delta-\langle\nabla w,\nabla\cdot\rangle$, it follows from integration by
	parts and \eqref{operador} that
	\begin{equation*}
	\int_{M}|\nabla h|^{2} e^{-w}dv=-\int_{M} h(\Delta_{w}h) e^{-w}dv=0.
	\end{equation*}
	Hence, $|\nabla h|=0$, which show that $h$ is constant.
	
\end{myproof}

\begin{myproof}{Theorem}{\ref{teorema2}} Using lemma \ref{lemma}, we have that 
\begin{equation}\label{11}
S_{B}+\frac{\lambda_{F}}{f^2}-\frac{2d}{f}\Delta f-d(d-1)\frac{|\nabla f|^2}{f^2}-\rho= \frac{\langle\nabla f, \nabla h\rangle}{f}.
\end{equation}

Since $S_{B}\geq \rho-\frac{\lambda_{F}}{f^2}$, we obtain by equation \eqref{11}, that
\begin{equation}\label{maximo principio}
\Delta f+\langle\nabla w,\nabla f\rangle=\frac{(S_{B}-\rho)f^2+\lambda_{F}}{2 fd}\geq0,
\end{equation}
where $w=\frac{h}{2d}+\ln f^{\frac{d-1}{2}}$

Now, consider $x_{0}$ the point where $f$ attains its maximum $f_{0}$, and define
\begin{equation*}
\Omega_{0}:=\{x\in B\hspace{0.2cm} ;\hspace{0.2cm} f(x)=f_{0}\}
\end{equation*} 
$\Omega_{0}$ is closed and nonempty since $x_{0}\in \Omega_{0}$. Let now $y\in \Omega_{0}$, then applying the maximum principle (see \cite{gilbarg2015elliptic} p. 35 ) to \eqref{maximo principio} we obtain that, $f(x)=f_{0}$ in a neighborhood of $y$ so that $\Omega_{0}$ is open. Connectedness of $B$ yields $\Omega_{0}=B$. Thus $f$ is constant.
\end{myproof}

\begin{myproof}{Proposition}{\ref{proposicao}}  Consider $\theta=\langle \nabla\log f, \nabla h\rangle$, then applying Lemma \ref{lemma}, we obtain
$\theta g_{B}=Hess (h)$. The result follow by lemma:
	
	\begin{lema}(\cite{pigola2011remarks}, Theorem 1)Let $(N^{n},\hat{g})$ be a complete manifold. Suppose that there exist a smooth function $h:N\rightarrow\mathbb{R}$ satisfying $Hess (h)=\theta\hat{g}$ for some constant $\theta\neq0$. Then $N^{n}$ is isometric to $\mathbb{R}^{n}$.
	\end{lema}
	
	\end{myproof}

\begin{myproof}{Theorem}{\ref{teorema invariancia geral}} The equivalence given by Lemma \ref{lemma} states that a necessary and sufficient condition to $\overline{M}^{n+d}$ be a gradient Yamabe soliton with potential function $\widetilde{h}$ is
\begin{equation}\label{equivalence1}
S_{B}+\frac{\lambda_{F}}{f^2}-\frac{2d}{f}\Delta f-d(d-1)\frac{|\nabla f|^2}{f^2}-\rho= \frac{\langle\nabla f, \nabla h\rangle}{f},
\end{equation}
and
\begin{equation}\label{equivalence2}
\left(S_{B}+\frac{\lambda_{F}}{f^2}-\frac{2d}{f}\Delta f-d(d-1)\frac{|\nabla f|^2}{f^2}-\rho\right)g_{B}= Hess(h).
\end{equation}
	
	We will use the above equations in combination with the invariant solutions technique to obtain equations \eqref{eq:09}, \eqref{eq1}, \eqref{eq2}, \eqref{eq:10} and \eqref{eq10}.
	
	First, for an arbitrary choice of a non zero vector $\alpha=(\alpha_{1},\dots,\alpha_{n})$, consider  $\xi:\mathbb{R}^{n}\rightarrow\mathbb{R}$ given by $\xi(x_{1},\dots,x_{n})=\alpha_{1}x_{1}+\dots+\alpha_{n}x_{n}$. Since we are assuming that $\varphi(\xi)$, $h(\xi)$ and $f(\xi)$ are functions of $\xi$, then we have
	\begin{equation}\label{invariante}
	\begin{aligned}
	\varphi_{,x_{i}}&=\varphi'\alpha_{i},\hspace{0.2cm} & f_{,x_{i}}&=f'\alpha_{i},\hspace{0.2cm} & h_{,x_{i}}&=h'\alpha_{i},\hspace{0.2cm} \\[10pt]
	\varphi_{,x_{i}x_{j}}&=\varphi''\alpha_{i}\alpha_{j},\hspace{0.2cm} & f_{,x_{i}x_{j}}&=f''\alpha_{i}\alpha_{j},\hspace{0.2cm} & 	h_{,x_{i}x_{j}}&=h''\alpha_{i}\alpha_{j}.\hspace{0.2cm}
	\end{aligned}
	\end{equation}

	It is well known that for the conformal metric  $g_{B}=\varphi^{-2}\delta$, the Ricci curvature is given by (\cite{besse2007einstein}):	
\begin{equation}\label{Ricci}Ric_{g_{B}}=\frac{1}{\varphi^{2}}\Big{\{}(n-2)\varphi Hess_{\delta}(\varphi)+[\varphi\Delta_{\delta}\varphi-(n-1)|\nabla_{\delta}\varphi|^{2}]g\Big{\}}.
\end{equation}
	And then combining \eqref{invariante} with \eqref{Ricci} we easily see that the scalar curvature on conformal metric
	is given by
	\begin{eqnarray}\label{escalar conforme}S_{B}&=&\sum_{k=1}^{n}\varphi^2\varepsilon_{k}(Ric_{g_{B}})_{kk}\nonumber\\
	&=&(n-1)(2\varphi\sum_{k=1}^{n}\varepsilon_{k}\varphi_{,x_{k}x_{k}}-n\sum_{k=1}^{n}\varepsilon_{k}\varphi_{,x_{k}}^{2})\nonumber\\
	&=&||\alpha||^{2}(n-1)(2\varphi\varphi''-n(\varphi')^{2}).
	\end{eqnarray}

	Now, in order to compute the $Hess (h)$ of $h$ relatively to $g_{B}$ we evoke the expression
	\begin{equation}
	(Hess(h))_{ij}=h_{,x_{i}x_{j}}-\sum_{k=1}^{n}\Gamma_{ij}^{k}h_{,x_{k}},
	\end{equation}
	where the Christoffel symbol $\Gamma_{ij}^{k}$ for distinct $i,j,k$ are given by
	\begin{equation}\bar{\Gamma}_{ij}^{k}=0,\ \bar{\Gamma}_{ij}^{i}=-\frac{\varphi_{,x_{j}}}{\varphi},\ \bar{\Gamma}_{ii}^{k}=\varepsilon_{i}\varepsilon_{k}\frac{\varphi_{,x_{k}}}{\varphi}\;\ \mbox{and}\;\ \bar{\Gamma}_{ii}^{i}=-\frac{\varphi_{,x_{i}}}{\varphi}.\nonumber
	\end{equation}
	Therefore,
		\begin{eqnarray}\label{hessian}(Hess(h))_{ij}&=&h_{,x_{i}x_{j}}+\varphi^{-1}(\varphi_{,x_{i}}h_{,x_{j}}+\varphi_{,x_{j}}h_{,x_{i}})-\delta_{ij}\varepsilon_{i}\sum_{k}\varepsilon_{k}\varphi^{-1}\varphi_{,x_{k}}h_{,x_{k}}\nonumber\\
	&=&\alpha_{i}\alpha_{j}h''+(2\alpha_{i}\alpha_{j}-\delta_{ij}\varepsilon_{i}||\alpha||^{2})\varphi^{-1}\varphi'h'.
	\end{eqnarray}
And the Laplacian $\Delta f=\sum_{k}\varphi^{2}\varepsilon_{k}(Hess(f))_{kk}$ of $f$ with respect to $g_{B}$ is
\begin{equation}
\Delta f=||\alpha||^{2}\varphi^{2}(f''-(n-2)\varphi^{-1}\varphi'f').
\end{equation}

On the other hand, the expression of $\langle\nabla f,\nabla h\rangle$ and $|\nabla f|^{2}$ on conformal metric $g_{B}$ are given by
\begin{equation}\label{gradientes}
\langle\nabla f,\nabla h\rangle=\varphi^{2}\sum_{k}\varepsilon_{k}f_{,x_{k}}h_{,x_{k}}=||\alpha||^{2}\varphi^{2}f'h',\quad |\nabla f|^{2}=\varphi^{2}\sum_{k}\varepsilon_{k}f_{,x_{k}}^{2}=||\alpha||^{2}\varphi^{2}(f')^{2}.
\end{equation}
Then substituting \eqref{escalar conforme} and \eqref{gradientes} into \eqref{equivalence1} we obtain \eqref{eq2}.
	
	Now, for $i\neq j$ we obtain by \eqref{equivalence2} and \eqref{hessian} that 
	\begin{equation*}
	\alpha_{i}\alpha_{j}\left(h''+2\frac{\varphi'h'}{\varphi}\right)=0.
	\end{equation*}
	
	If there exist $i,j$, $i\neq j$ such that $\alpha_{i}\alpha_{j}\neq
	0$, then we get
	\begin{equation*}h''+2\frac{h'\varphi'}{\varphi}=0
	\end{equation*}
which is equation \eqref{eq:09}. And for $i=j$, substituting \eqref{escalar conforme}, \eqref{hessian} and \eqref{gradientes} into \eqref{equivalence2} we obtain \eqref{eq1}. 

Now, we need to consider the case $\alpha_{k_{0}}=1$, $\alpha_{k}=0$ for $k\neq k_{0}$. In this case, substituting \eqref{escalar conforme}, \eqref{hessian} and \eqref{gradientes} into \eqref{equivalence2} we obtain
\begin{equation*}
\begin{split}
\Big{[}\varepsilon_{k_{0}}(n-1)(2\varphi\varphi''-n(\varphi')^{2})+\frac{\lambda_{F}}{f^{2}}-\frac{2d}{f}\varepsilon_{k_{0}}(\varphi^{2}f''-(n-2)\varphi\varphi'f')+\\
-\varepsilon_{k_{0}}\frac{d(d-1)}{f^2}\varphi^{2}(f')^{2}-\rho\Big{]}\frac{\varepsilon_{i}}{\varphi^{2}}=-\varepsilon_{i}\varepsilon_{k_{0}}\frac{\varphi'}{\varphi}h'
\end{split}
\end{equation*}
when $i\neq k_{0}$, that is, $\alpha_{i}=0$

and
\begin{equation*}
\begin{split}
\Big{[}\varepsilon_{k_{0}}(n-1)(2\varphi\varphi''-n(\varphi')^{2})+\frac{\lambda_{F}}{f^{2}}-\frac{2d}{f}\varepsilon_{k_{0}}(\varphi^{2}f''-(n-2)\varphi\varphi'f')+\\
-\varepsilon_{k_{0}}\frac{d(d-1)}{f^2}\varphi^{2}(f')^{2}-\rho\Big{]}\frac{\varepsilon_{k_{0}}}{\varphi^{2}}=h''+\frac{\varphi'}{\varphi}h'
\end{split}
\end{equation*}
for $i=k_{0}$, that is, $\alpha_{k_{0}}=1$.

However, this equations are equivalent to equations \eqref{eq:09}
and $\eqref{eq1}$. The Lightlike case follow by taking $||\alpha||^{2}=0$ into \eqref{eq1} and \eqref{eq2}. This completes the demonstration.

\end{myproof}

\begin{myproof}{Theorem}{\ref{eq5} and \ref{eq5'}}
	Since $\rho=0$ and $h'\neq0$ we have by equation \eqref{eq1} and \eqref{eq2}
	of Theorem \ref{teorema invariancia geral} that
	\begin{equation*}\frac{\varphi'}{\varphi}=-\frac{f'}{f}.
	\end{equation*}
	Integrating this equation we have
	\begin{equation}\label{eq14}f(\xi)=\frac{k_{2}}{\varphi(\xi)},
	\end{equation}
	for some $k_{2}\in\mathbb{R}\setminus\{0\}$, which is equation \eqref{eqa} and \eqref{eqa'}.
	
	Integrating the equation \eqref{eq:09}, we have that
	\begin{equation}\label{eq19}h'(\xi)=\frac{k_{1}}{\varphi^{2}(\xi)},
	\end{equation}
	for some $k_{1}\neq0$, and
	$$h(\xi)=k_{1}\int\frac{1}{\varphi^{2}(\xi)}d\xi,$$
	which is equation \eqref{eqb} and \eqref{eqb'}.
	
	Substituting equation \eqref{eq19} into \eqref{eq1} and considering \eqref{eq14}, we obtain the follow differential equation
	\begin{equation*}\varphi^2\varphi''-\frac{(n+d)}{2}\varphi(\varphi')^{2}+\frac{k_{1}}{2(n+d-1)}\varphi'=-\frac{\varphi^{3}\lambda_{F}}{k_{2}^{2}||\alpha||^{2}}.
	\end{equation*}
	Then since, $n+d=6$, we obtain
	\begin{equation}\label{eq16}\varphi^2\varphi''-3\varphi(\varphi')^{2}+\frac{k_{1}}{10}\varphi'=-\frac{\varphi^{3}\lambda_{F}}{k_{2}^{2}||\alpha||^{2}}.
	\end{equation}

	If $\lambda_{F}=0$, considering $\varphi(\xi)^{-2}=v(\xi)$, we obtain by \eqref{eq16} the following equivalent condition
	\begin{equation}\label{eq177}v''+\frac{k_{1}}{10}v'v=0.
	\end{equation}
	Integrating equation \eqref{eq177} we get
	\begin{equation*}v'+\frac{k_{1}}{20}v^{2}=k_{2},\quad k_{2}=\text{constant}.
	\end{equation*}
	This implies that
	\begin{equation}-\int\frac{1}{\frac{k_{1}}{20}v^{2}-k_{2}}dv=\xi+k_{3},\quad k_{3}=\text{constant}.\nonumber
	\end{equation}
	Therefore, it follows from $\varphi(\xi)^{-2}=v(\xi)$ that
	\begin{equation}40\int\frac{\varphi d\varphi}{k_{1}-20k_{2}\varphi^{4}}=\xi+k_{3},\nonumber
	\end{equation}
	which is equation \eqref{eqc'}.

	For $\lambda_{F}\neq0$, consider the change $u(\varphi)=\varphi^{-3}\varphi'$, then \eqref{eq16} is equivalent to the following separable variables differential equation
	\begin{equation*}\frac{du}{d\varphi}=\frac{-q-pu}{\varphi^{5}u},
	\end{equation*}
	where $p=\frac{k_{1}}{10}$, $q=\frac{\lambda_{F}}{k_{2}^{2}||\alpha||^{2}}$, and which solution is
given by	
	\begin{equation*}
	u(\varphi)=\frac{q}{p}\big(W\big(k_{3}e^{\frac{p^{2}}{4q\varphi^{4}}-1}\big)+1\big),\quad k_{3}=\text{constant}\neq 0,
	\end{equation*}
	where  $W$ is the product log function.
	
	Substituting back for $u(\varphi)=\varphi^{-3}\varphi'$, we obtain
	\begin{equation*}
	\frac{pd\varphi}{q\varphi^3W\big(k_{3}e^{\frac{p^{2}}{4q\varphi^4}-1}\big)+q\varphi^3}=dt,
	\end{equation*}
	which provide equation \eqref{eqc}. The converse is a straightforward
	computation. This concludes the proof of Theorems \ref{eq5} and \ref{eq5'}.

\end{myproof}

\begin{myproof}{Theorem}{\ref{eq7}}Since $h'=0$ and $\lambda_{F}=\rho=0$ we have by equation \eqref{eq1} and \eqref{eq2} of Theorem \ref{teorema invariancia geral} that
	\begin{equation}(n-1)(2\varphi\varphi''-n(\varphi')^{2})-2\frac{d}{f}(\varphi^{2}f''-(n-2)\varphi\varphi'f')-\frac{d(d-1)}{f^2}\varphi^{2}(f')^{2}=0,\nonumber
	\end{equation}
	which is equivalent to
	\begin{equation}\left(\frac{f'}{f}-\frac{(n-2)}{(d+1)}\frac{\varphi'}{\varphi}\right)^{2}+\frac{2}{d+1}\left(\frac{f'}{f}-\frac{(n-2)}{(d+1)}\frac{\varphi'}{\varphi}\right)^{'}+\frac{(n+d-1)}{d(d+1)^{2}}\left(n\left(\frac{\varphi'}{\varphi}\right)^{2}-2\frac{\varphi''}{\varphi}\right)=0.\nonumber
	\end{equation}
	Consider
	$z=\frac{f'}{f}-\frac{(n-2)}{(d+1)}\frac{\varphi'}{\varphi}$, then
	\begin{equation}\label{eq20}z^{2}+\frac{2}{d+1}z'+\frac{(n+d-1)}{d(d+1)^{2}}\left(n\left(\frac{\varphi'}{\varphi}\right)^{2}-2\frac{\varphi''}{\varphi}\right)=0.
	\end{equation}
	
	Recall that the Ricatti differential equation is a differential
	equation of the form
	\begin{equation}\label{eq21}z(\xi)'=f_{2}(\xi)z(\xi)^{2}+f_{1}(\xi)z(\xi)+f_{0}(\xi),
	\end{equation}
	where $f_{0}$, $f_{1}$ and $f_{2}$ are smooth functions on $\mathbb{R}$, and by
	Picard theorem, given a particular solution $z_{0}$ of \eqref{eq21} we have that the general solution of Riccati equation is given by
\begin{equation*}
	z(\xi)=z_{0}(\xi)+\Phi(\xi)\left[C-\int\Psi(\xi)f_{2}(\xi)d\xi\right]^{-1},
\end{equation*}
where
\begin{equation*}
\Psi(\xi)=\exp\Big{\{}\int[2f_{2}(\xi)z_{0}(\xi)+f_{1}(\xi)]d\xi\Big{\}},\quad C=\text{constant}.
\end{equation*}

	Observe that \eqref{eq20} is a Ricatti differential equation with
	\begin{equation*}f_{1}(\xi)=0,\hspace{0,2cm} f_{2}(\xi)=-\frac{d+1}{2}\hspace{0,2cm} \text{and}\hspace{0,2cm}
	f_{0}(\xi)=-\frac{(n+d-1)}{2d(d+1)}\left(n\left(\frac{\varphi'}{\varphi}\right)^{2}-2\frac{\varphi''}{\varphi}\right).
	\end{equation*}
	Then we obtain
	$$\frac{f'(\xi)}{f(\xi)}=\frac{(n-2)}{(d+1)}\frac{\varphi'(\xi)}{\varphi(\xi)}+z_{p}(\xi)+\frac{e^{-(d+1)\int z_{p}(\xi)d\xi}}{\frac{d+1}{2}\int e^{-(d+1)\int z_{p}(\xi)d\xi}d\xi+C}.$$
	And thus
	$$f(\xi)=\varphi^{\frac{n-2}{d+1}}e^{\int z_{p}d\xi}\left(\int e^{-(d+1)\int z_{p}d\xi}d\xi+\frac{2}{d+1}C\right)^{\frac{2}{d+1}},$$
	where $z_{p}(\xi)$ is a particular solution of \eqref{eq20}. This
	expression is equation \eqref{eq22}.
	
	Now, since $h'=0$, we have that $h(\xi)=constant$, which is equation
	\eqref{eq23}. Then we prove the necessary condition. A direct calculation shows us the converse implication. This
	concludes the proof of Theorem.

\end{myproof}

\begin{myproof}{Theorem}{\ref{eq6}} It follows immediately from \eqref{eq:10} and \eqref{eq10}.

\end{myproof}

\begin{myproof}{completeness of example}{\ref{exemplo completo}}
	
Let $(\mathbb{R}^{n},g)$ be the standard semi-Euclidean space where $g=-dx_{1}^{2}+\sum_{i=2}^{n}dx_{i}^{2}$. Take $k\in\mathbb{R}\setminus\{0\}$, and consider the functions \begin{equation}\varphi(\xi)=e^{k\xi}, \hspace{0,3cm}f(\xi)=e^{k\xi},\hspace{0,3cm}h(\xi)=-\frac{k_{1}e^{-2k\xi}}{2k},\quad k_{1}\neq0.\nonumber
	\end{equation}

Call $\hat{g}:=\varphi^{-2}g=e^{-2k\xi}g$, then the gradient $\nabla_{\hat{g}}f$ is given by
	\begin{equation*}
	\nabla_{\hat{g}}f=\sum_{r,s=1}^{n}\hat{g}^{rs}f_{,x_{s}}\partial_{r}=\sum_{r,s=1}^{n}\varphi^{2}\varepsilon_{r}\delta_{rs}f'\alpha_{s}\partial_{s}=\sum_{s=1}^{n}k\varepsilon_{s}\alpha_{s}e^{3k\xi}\partial_{s}.
	\end{equation*}
Since $\alpha_{1}=\alpha_{2}=1$, $\alpha_{i}=0$, for $i\geq3$, and $\varepsilon_{1}=-1$, $\varepsilon_{i}=1$, for $i\geq2$, we obtain
	\begin{equation*}
\nabla_{\hat{g}}f=\big{(}-ke^{3k\xi},ke^{3k\xi}, 0,\dots, 0\big{)}.
\end{equation*}
Then, considering $\gamma_{B}(s)=(y_{1}(s),\dots, y_{n}(s))$ and $\gamma_{F}(s)=(y_{n+1},\dots, y_{n+p}(s))$ in Proposition \ref{geodesicas produto torcido}, we have
\begin{equation*}
\left\{\begin{array}{r@{\mskip\thickmuskip}l}
	y_{1}''(s)&= -k\big{[}y_{n+1}'(s)^{2}+\dots+y_{n+p}'(s)^{2}\big{]}e^{4k(y_{1}(s)+y_{2}(s))}, \hfill \text{(I)}\\ [10pt]
	y_{2}''(s)&= k\big{[}y_{n+1}'(s)^{2}+\dots+y_{n+p}'(s)^{2}\big{]}e^{4(y_{1}(s)+y_{2}(s))}, \hfill \text{(II)}\\[10pt]
	y_{r}''(s)&= 0,\quad \text{for}\quad r\in\{3,\dots,n\},\hfill \text{(III)}\\[10pt]
	y_{n+l}''(s)&=-2k[y_{1}'(s)+y_{2}'(s)]y_{n+l}'(s),\quad \text{for}\quad l\in\{1,\dots,d\}.\qquad \text{(IV)}
\end{array} \right.
\end{equation*}

The sum of differential equation \text{(I)} and \text{(II)} gives $y_{1}''(s)+y_{2}''(s)=0$, then by integration 
\begin{equation}\label{diferential eq1}
y_{1}'(s)+y_{2}'(s)=c_{1}, \quad y_{1}(s)+y_{1}(s)=c_{1}s+c_{2}, \quad c_{1}, c_{2}\in \mathbb{R}.
\end{equation}

Substituting \eqref{diferential eq1} into \text{(IV)}, we obtain the second order linear ordinary differential equation
\begin{equation}\label{second order}
y_{n+l}''(s)+2kc_{1}y_{n+l}'(s)=0 \quad \text{for each} \quad l\in\{1,\dots,d\},
\end{equation}
whose general solutions is
\begin{equation*}
y_{n+l}(s) = \begin{cases} c_{3,l}+c_{4,l}s &\text{if } c_{1}=0 \\[10pt]
	c_{3,l}+c_{4,l}e^{-2kc_{1}s} & \text{if } c_{1}\neq0 
\end{cases}
\end{equation*}
where $c_{3,l},c_{4,l}\in\mathbb{R}$. This shows that for each $l\in\{1,\dots,d\}$, the functions $y_{n+l}(s)$ are defined on the entire real line $\mathbb{R}$. Since the solutions of \text{(III)} are given by $y_{r}(s)=c_{5,r}+c_{6,r}s$, for $c_{5,r}, c_{6,r}\in\mathbb{R}$, whose domain is $\mathbb{R}$, it is only necessary to prove that the solutions of \text{(I)} and \text{(II)} are also defined in $\mathbb{R}$.

Integrating \eqref{second order} and replacing its result into \text{(I)}, we have
\begin{equation*}
y_{1}''(s)=-k[c_{7,1}^{2}+c_{7,2}^{2}+\dots+c_{7,r}^{2}]e^{-4kc_{1}s}e^{4k(y_{1}(s)+y_{2}(s))}
\end{equation*}
where $c_{7,1}, c_{7,2},\dots,c_{7,r}\in\mathbb{R}$.

Now, by \eqref{diferential eq1} we obtain that

\begin{eqnarray}
y_{1}''(s)&=&-k[c_{7,1}^{2}+c_{7,2}^{2}+\dots+c_{7,p}^{2}]e^{-4kc_{1}s}e^{4k(c_{1}s+c_{2})}\nonumber\\
&=&-k[c_{7,1}^{2}+c_{7,2}^{2}+\dots+c_{7,p}^{2}]e^{4kc_{2}}\nonumber\\
&=& c_{8,1}\in\mathbb{R}.
\end{eqnarray}
Then $y_{1}(s)=\frac{c_{8,1}}{2}s^{2}+c_{9}s+c_{10}$, whose domain is $\mathbb{R}$. Except for the signal, the same occurs for $y_{2}(s)$. Thus, all the geodesics $\gamma=(\gamma_{B},\gamma_{F})$ are
defined for the entire real line, which means that $(\mathbb{R}^{n},\varphi^{-2}g)\times_{f}(\mathbb{R}^{d},g_{0})$ is geodesically complete.
\end{myproof}


\begin{thebibliography}{X}



\bibitem{yano1957theory} Yano, Kentaro L. \textit{The theory of Lie derivatives and its applications.} North-Holland, (1957)

\bibitem{hamilton1988ricci} Hamilton, Richard S. \textit{The Ricci flow on surfaces, Mathematics and general relativity (Santa Cruz, CA, 1986), 237--262}. Contemp. Math. \textbf{71}, 301--307,
(1988)

\bibitem{brozos2016local} Brozos-V{\'a}zquez, Miguel and Calvi{\~n}o-Louzao, Esteban and Garc{\'\i}a-R{\'\i}o, Eduardo and V{\'a}zquez-Lorenzo, Ram{\'o}n. \textit{Local structure of self-dual gradient Yamabe solitons.} Geometry, Algebra and Applications: From Mechanics to Cryptography. Springer, 25--35, (2016)

\bibitem{calvino2012three} Calvi{\~n}o-Louzao, E and Seoane-Bascoy, J and V{\'a}zquez-Abal, ME and V{\'a}zquez-Lorenzo, R. \textit{Three-dimensional homogeneous Lorentzian Yamabe solitons.}
Abhandlungen aus dem Mathematischen Seminar der Universit{\"a}t Hamburg, vol 82(2), 193--203, (2012)

\bibitem{o1983semi} O'neill, Barrett. \textit{Semi-Riemannian geometry with applications to relativity.} vol 103, Academic press, (1983)

\bibitem{al2012warped} Al-Solamy, Falleh R and Khan, Meraj Ali. \textit{Warped product submanifolds of Riemannian product manifolds.} Abstract and Applied Analysis, \textbf{2012}, (2012)

\bibitem{blaga2017warped} Blaga, Adara M. \textit{On warped product gradient $eta $-Ricci solitons.} arXiv preprint arXiv:1705.04092, (2017)

\bibitem{feitosa2017construction} Feitosa, FES and Freitas Filho, AA and Gomes, JNV. \textit{On the construction of gradient Ricci soliton warped product.} Nonlinear Analysis, \textbf{161},
30--43, (2017)

\bibitem{kim2013warped} Kim, Byung Hak and Lee, Sang Deok and Choi, Jin Hyuk and Lee, Young Ok. \textit{On warped product spaces with a certain Ricci condition.} Bulletin of the Korean Mathematical Society, \textbf{50(5)}, 1683--1691, (2013)

\bibitem{de2017gradient} de Sousa, M{\'a}rcio Lemes and Pina, Romildo. \textit{Gradient Ricci solitons with structure of warped product.} Results in Mathematics, \textbf{71(3-4)}, 825--840,
(2017)

\bibitem{leandro2017invariant} Leandro, Benedito and Pina, Romildo. \textit{Invariant solutions for the static vacuum equation.} Journal of Mathematical Physics, \textbf{58(7)}, 072502, 10 pp. (2017)

\bibitem{olver2000applications} Olver, Peter J. \textit{Applications of Lie groups to differential equations.} Springer, \textbf{107}, (2000)

\bibitem{he2011gradient} He, Chenxu. \textit{Gradient Yamabe solitons on warped products.} arXiv preprint arXiv:1109.2343, (2011)

\bibitem{neto2018gradient} Neto, Benedito Leandro and Tenenblat, Keti. \textit{On gradient Yamabe solitons conformal to a pseudo-Euclidian space.} Journal of Geometry and Physics, \textbf{123}, 284--291, (2018)

\bibitem{barbosa2014gradient} Barbosa, Ezequiel and Pina, Romildo and Tenenblat, Keti. \textit{On gradient Ricci solitons conformal to a pseudo-Euclidean space.} Israel Journal of Mathematics, \textbf{200(1)}, 213--224, (2014)

\bibitem{barboza2018invariant} Barboza, Marcelo and Leandro, Benedito and Pina, Romildo. \textit{Invariant solutions for the Einstein field equation.} Journal of Mathematical Physics, \textit{59(6)}, 062501, 9 pp. (2018)

\bibitem{fernandez2001curvature} Fern{\'a}ndez-L{\'o}pez, M and Garcia-Rio, E and Kupeli, DN and {\"U}nal, B. \textit{A curvature condition for a twisted product to be a warped product.} manuscripta mathematica, \textbf{106(2)}, 213--217, (2001)

\bibitem{bishop1969manifolds} Bishop, Richard L and O'Neill, Barrett. \textit{Manifolds of negative curvature.} Transactions of the American Mathematical Society, \textbf{145}, 1--49, (1969)

\bibitem{ganchev2000riemannian} Ganchev, G and Mihova, V. \textit{Riemannian manifolds of quasi-constant sectional curvatures.} Journal fur die Reine und Angewandte Mathematik, 119--142, (2000)

\bibitem{dobarro1987scalar} Dobarro, F and Lami Dozo, E. \textit{Scalar curvature and warped products of Riemann manifolds.} Transactions of the American Mathematical Society, \textbf{303(1)}, 161--168, (1987)

\bibitem{besse2007einstein} Besse, Arthur L. \textit{Einstein manifolds.} Springer, (2007)

\bibitem{lopez2016geometry} L{\'o}pez, Marco Castrill{\'o}n and Encinas, Luis Hern{\'a}ndez and Gadea, Pedro Mart{\'\i}nez and Mar{\'\i}a, Ma Eugenia Rosado. \textit{Geometry, Algebra and Applications: From Mechanics to Cryptography.} Springer, (2016)

\bibitem{batat2014ricci} Batat, Wafaa and Onda, Kensuke. \textit{Ricci and Yamabe solitons on second-order symmetric, and plane wave 4-dimensional Lorentzian manifolds.} Journal of Geometry,
\textbf{105(3)}, 561--575, (2014)

\bibitem{lee2017warped} Lee, Sang Deok and Kim, Byung Hak and Choi, Jin Hyuk. \textit{Warped product spaces with Ricci conditions.} Turkish Journal of Mathematics, \textbf{41(6)}, 1365--1375, (2017)

\bibitem{barbosa2013conformal} Barbosa, Ezequiel and Ribeiro, Ernani. \textit{On conformal solutions of the Yamabe flow.} Archiv der Mathematik, \textbf{101(1)}, 79--89, (2013)

\bibitem{cao2012structure} Cao, Huai-Dong and Sun, Xiaofeng and Zhang, Yingying. \textit{On the structure of gradient Yamabe solitons.} Mathematical Research Letters, \textbf{19(4)}, 767--774, (2012)

\bibitem{daskalopoulos2013classification} Daskalopoulos, Panagiota and Sesum, Natasa. \textit{The classification of locally conformally flat Yamabe solitons.} Advances in Mathematics, \textbf{240}, 346--369, (2013)

\bibitem{hsu2012note} Hsu, Shu-Yu. \textit{A note on compact gradient Yamabe solitons.} Journal of Mathematical Analysis and Applications, \textbf{388(2)}, 725--726, (2012)

\bibitem{bang2017differential} Bang-yen, Chen. \textit{Differential geometry of warped product manifolds and submanifolds.} World Scientific, (2017)

\bibitem{pigola2011remarks} Pigola, Stefano and Rimoldi, Michele and Setti, Alberto G. \textit{Remarks on non-compact gradient Ricci solitons.} Mathematische Zeitschrift, \textbf{268(3-4)}, 777--790, (2011)

\bibitem{gilbarg2015elliptic} Gilbarg, David and Trudinger, Neil S. \textit{Elliptic partial differential equations of second order.} Springer, (2015)


\end{thebibliography}
\end{document}